\documentclass[11pt,a4paper]{article}
\usepackage[centertags]{amsmath}
\usepackage{amsfonts}
\usepackage{amssymb}
\usepackage{amsthm}
\usepackage{newlfont}

\newlength{\defbaselineskip}
\theoremstyle{plain}
\newtheorem{thm}{Theorem}[section]

\newtheorem{defn}{Definition}[section]

\newtheorem{rem}{Remark}[section]
\numberwithin{equation}{section}

\newcommand{\A}{\alpha}
\newcommand{\B}{\beta}

\newcommand{\MM}{\widetilde{M}}

\newcommand{\G}{\widetilde{g}}

\newcommand{\n}{\nabla}

\newcommand{\PB}{\overline{P}}
\newcommand{\PP}{\widetilde{P}}

\title{Induced structures on product of spheres}
\author{Cristina-Elena Hre\c{t}canu}
\date{}

\usefont{T1}{cmss}{m}{ol}
\begin{document}

\maketitle

\begin{abstract}
The purpose of this paper is to give an effective construction for
some induced structures on spheres or product of spheres of
codimension $1$, $2$ and $3$, respectively, in Euclidean space
endowed with an almost product structure.

\end{abstract}

\section{Introduction}
\normalfont Let $(\MM,\G)$ be a Riemannian manifold equipped with
a Riemannian metric tensor $\G$ and a (1,1) tensor field $\PP$
such that
\[\PP^{2} = \epsilon I, \leqno(1.1)\]
where $I$ is the identity on $\MM$ and $\epsilon =\pm 1$.  We
suppose that $\G$ and $\PP$ are compatible in the sense that we
have
\[\G(\PP U, \PP V)=\G(U, V), \leqno(1.2) \]
which is equivalent to
\[ \G(\PP U, V)=\epsilon\G(U, \PP V), \leqno(1.3)\]
for each $U, V \in \chi(\MM)$, where $\chi(\MM)$ is the Lie
algebra of the vector fields on $\MM$.

For $\epsilon =1$ we obtain that $\PP$ is an almost product
structure and the Riemannian manifold $(\MM,\G,\PP)$, with the
compatibility relation (1.2), becomes an almost product Riemannian
manifold.

Let $M$ be a $n$-dimensional submanifold of codimension $r$ $(n,r
\in \mathbb{N}^{*})$ in a Riemannian manifold $(\MM,\G,\PP)$ which
satisfied the relations (1.1) and (1.2).

We note the tangent space of M in a point $x \in M$ by $T_{x}(M)$
and the normal space of M in x by $T_{x}^{\bot}(M)$. For $\A \in
\{1,...,r\}$, let $(N_{1},...,N_{r}):=(N_{\A})$ be an orthonormal
basis in $T_{x}^{\bot}(M)$, for every $x \in M$.

In the same manner like in \cite{Adati2} and \cite{eu7}, we
construct the $(a,\epsilon)f$ Riemannian structure, which is a
structure $(P,g,\epsilon \xi_{\A},u_{\A},(a_{\A\B})_{r})$, induced
on a submanifold M in a Riemannian manifold $(\MM,\G,\PP)$ ($\PP$
and $\G$ verify the relation (1.1) and (1.2)), where P is an
(1,1)-tensor field on M, g is the induced metric on submanifold M,
$\xi_{\A}$ are tangent vector fields on $M$, $u_{\A}$ are 1-forms
on $M$ and $(a_{\A\B})_{r}$ is a $r \times r$ matrix of real
functions on $M$.

The $(a,\epsilon)f$ Riemannian structure induced on M by the $\PP$
structure on $(\MM,\G)$ (which verifies the relations (1.1) and
(1.2)) is a generalization of the almost r-paracontact Riemannian
structure(\cite{Bucki2}). The case of the $(a,-1) f$ Riemannian
structure was studied by K. Yano and M. Okumura (in \cite{Yano1}
and \cite{Yano2}).

%>>>>>>>>>>>>>>>>>>>>>>>>>>>>>>>>>>>>>>>>>>>>>>>>>>>>>>>>>>>>>>>>
\section {$(a,\epsilon)f$ Riemannian structure}

Let M be a $n$-dimensional submanifold of codimension $r$ ($n,r
\in \mathbb{N}^{*}$) in a Riemannian manifold $(\MM, \G)$,
equipped by an (1,1)-tensor field $\PP$, such that $\G$ and $\PP$
verify the conditions (1.1) and (1.2). The decomposition of vector
fields $\PP X$ and $\PP N_{\A}$, respectively, in the tangential
and normal components on the submanifold $M$ in $\MM$ has the
form:
\[ \PP X = P X + \sum_{\A=1}^{r}u_{\A}(X)N_{\A},  \leqno(2.1)\]
for every $X \in \chi(M)$ and
\[ \PP N_{\A}= \epsilon \xi_{\A} + \sum_{\B=1}^{r}a_{\A \B}
N_{\B}, \quad(\epsilon = \pm 1)\leqno(2.2)\] for every $\A \in
\{1,...,r\}$.

\normalfont We use the notation $a:=(a_{\A\B})_{r}$ for the matrix
which result from the relation (2.2).

\defn   An $(a,\epsilon)f$
Riemannian structure on a n-dimensional submanifold M of
codimension r ($n,r \in \mathbb{N}^{*}$) in a Riemannian manifold
$(\MM,\G)$ is a structure $(P,g,\epsilon \xi_{\A},u_{\A},(a_{\A
\B})_{r})$ induced on M by the structure $(\PP,\G)$ from $\MM$
(where $\G$ and $\PP$ verify the conditions (1.1) and (1.2)).

\thm (\cite{eu7}) Let  M be a $n$-dimensional submanifold of
codimension $r$ ($n,r \in \mathbb{N}^{*}$) in a Riemannian
manifold $(\MM, \G)$, equipped by an (1,1)-tensor field $\PP$,
such that $\G$ and $\PP$ verify the conditions (1.1) and (1.2).
The structure $(\PP,\G)$  induces on the submanifold M a
$(a,\epsilon)f$ Riemannian structure, which verifies the following
properties:

\[(2.3)\begin{cases}
(i)\quad P^{2}X = \epsilon (X
-\sum_{\A=1}^{r}u_{\A}(X)\xi_{\A}),\\
(ii)\: u_{\A}(PX)= -\sum_{\B=1}^{r}a_{\B\A}u_{\B}(X),\\
(iii)\:\: a_{\A\B} = \epsilon a_{\B \A}, \\
(iv)\: u_{\A}(\xi_{\B})=
    \delta_{\A \B}-\epsilon \sum_{\gamma=1}^{r}a_{\A
    \gamma}a_{\gamma\B}, \\
(v)\quad P\xi_{\A}=-\sum_{\B=1}^{r}a_{\A \B}\xi_{\B}, \\
(vi)\quad u_{\A}(X)=g(X,\xi_{\A}),\\
(vii)\:\: g(PX,Y)=\epsilon g(X,PY),\\
(viii)\: g(PX,PY)= g(X,Y)-\\
\quad \quad -\sum_{\A=1}^{r}u_{\A}(X)u_{\A}(Y)
\end{cases}\]
 for any $X,Y \in \chi(M)$ and $\A,\B \in \{1,...,r\}$.

\rem \normalfont From $(2.3)(i),(v),(vi)$ we obtain
\[ \PP^{3}X=\epsilon P X +
\sum_{\A,\B=1}^{r}a_{\A\B}g(X,\xi_{\B})\xi_{\A} \leqno(2.4) \] for
every $X \in \chi(M)$.

Let $\widetilde{\n}$ and $\n$ be the Levi-Civita connections
defined on M and $\MM$ respectively, with respect to $\G$ and g
respectively.

\defn If $(\MM,\G,\PP)$ is an almost product Riemannian manifold
such that $\widetilde{\n} \PP =0$, then we say that $\MM$ is a
locally product Riemannian manifold.

\normalfont The Nijenhuis torsion tensor field of P has the form
$N_{P}(X,Y)=[PX,PY]+P^{2}[X,Y]-P[PX,Y]-P[X,PY]$ for any $X,Y \in
\chi(M)$. As in the case of an almost paracontact structure
(\cite{IMihai}), one can defined the normal
$(P,g,u_{\A},\xi_{\A},(a_{\A\B})_{r})$ structure on M.

\defn If we have the equality
\[ N_{P}(X,Y)-2\sum_{\A}du_{\A}(X,Y)\xi_{\A}=0, \leqno(2.5)\] for any $X,Y\in
\chi(M)$, then the $(a,1)f$ induced structure on submanifold M in
a Riemannian almost product manifold $(\MM,\G,\PP)$ is said to be
normal.

\thm (\cite{eu7}) Let M be an n-dimensional submanifold of
codimension r in a locally product Riemannian manifold
$(\widetilde{M},\G,\widetilde{P})$ and we suppose that the normal
connection $\n^{\bot}$ vanishes identically. If the determinant of
the matrix ($I_{r}-a^{2}$) does not vanish in any point $x \in M$,
then the $(a,1)f$ induced structure on M is normal if and only if
the (1,1) tensor field P commutes with the Weingarten operators
$A_{\A}$ (for every $\A \in \{1,...,r\}$).

\rem \normalfont If the submanifold M in $\MM$ is totally
umbilical, then the commutativity between the (1,1) tensor field P
and the Weingarten operators $A_{\A}$ (for any $\A \in
\{1,...,r\}$) has taken place. Therefore, if the determinant of
the matrix ($I_{r}-a^{2}$) does not vanish in any point $x \in M$
then the $(a,1)f$ induced structure on a totaly umbilical
submanifold is normal.

%>>>>>>>>>>>>>>>>>>>>>>>>>>>>>>>>>>>>>>>>>>>>>>>>>>>>>>>>>>>>>>>>>>>>>>>

\section{\bf Examples of $(a,1)f$ Riemannian structures
on spheres or product of spheres in Euclidean space }

 \normalfont

We suppose that the ambient space is $E^{2p+q}$ ($p,q \in
\mathbb{N}^{*}$) and for any point of $E^{2p+q}$ we have its
coordinates:
\[(x^{1},...,x^{p},y^{1},...,y^{p},z^{1},...,z^{q}):=(x^{i},y^{i},z^{j})\]
where $i\in \{1,...,p\}$ and $j \in \{1,...,q\}$. The tangent
space $T_{x}(E^{2p+q})$ is isomorphic with $E^{2p+q}$.

Let $\PP:E^{2p+q}\rightarrow E^{2p+q}$ an almost product structure
on $E^{2p+q}$ such that
\[\PP(x^{1},...,x^{p},y^{1},...,y^{p},z^{1},...,z^{q})=\leqno(3.1)\]
\[=(y^{1},...,y^{p},x^{1},...,x^{p},\varepsilon z^{1},...,\varepsilon z^{q})\] where $\varepsilon = \pm 1$.
Thus, $(\PP, <\:>)$ is an almost product Riemannian structure on
$E^{2p+q}$.

\textbf{ Example 1.} In this example, we construct an
$(a,1)f$-structure on the sphere $S^{2p+q-1}(R)\hookrightarrow
E^{2p+q}$. The equation of sphere $S^{2p+q-1}(R)$ is
\[\sum_{i=1}^{p}(x^{i})^{2}+\sum_{i=1}^{p}(y^{i})^{2}+
\sum_{j=1}^{q}(z^{j})^{2}=R^{2} \leqno(3.2)\] where R is its
radius and $(x^{1},...,x^{p},y^{1},...,y^{p},z^{1},...,z^{q})$ are
the coordinates of any point of $S^{2p+q-1}(R)$. We use the
following notations
\[\sum_{i=1}^{p}(x^{i})^{2}=r_{1}^{2},\:
\sum_{i=1}^{p}(y^{i})^{2}=r_{2}^{2},\:
\sum_{j=1}^{q}(z^{j})^{2}=r_{3}^{2},\leqno(3.3)\] and
$r_{1}^{2}+r_{2}^{2}=r^{2}$. Thus, we have
$r^{2}+r_{3}^{2}=R^{2}$. We remark that
\[N_{1}:= \frac{1}{R}(x^{i},y^{i},z^{j}),\:\: i \in {1,...,p},j \in \{1,...,q\} \leqno(3.4)\]
is a unit normal vector field on sphere $S^{2p+q-1}(R)$ and
\[\PP N_{1}=\frac{1}{R}(y^{i},x^{i},\varepsilon z^{j}) \leqno(3.5)\]

For a tangent vector field $X$ on $S^{2p+q-1}(R)$ we use the
following notation
\[X=(X^{1},...,X^{p},Y^{1},...,Y^{p},Z^{1},...,Z^{q}):=(X^{i},Y^{i},Z^{j})
\leqno(3.6)\] Hence we have
\[\sum_{i=1}^{p}x^{i}X^{i}+\sum_{i=1}^{p}y^{i}Y^{i}+\sum_{j=1}^{q}z^{j}Z^{j}=0 \leqno(3.7)\]
If we decompose $\PP N_{1}$ in the tangential and normal
components at sphere $S^{2p+q-1}(R)$, we obtain
\[ \PP N_{1}= \overline{\xi}+a_{11}N_{1} \leqno(3.8) \]
where $a_{11}= < \PP N_{1},N_{1}>$. If we use the following
notation
\[\sigma = \sum_{i=1}^{p}x^{i}y^{i} \leqno(3.9)\]
we obtain
\[a_{11}=\frac{2 \sigma + \varepsilon r_{3}^{2}}{R^{2}}. \leqno(3.10)\]
From (3.8) we obtain the tangential component of $\PP N_{1}$ at
sphere $S^{2p+q-1}(R)$:
\[\overline{\xi}=\frac{1}{R}(y^{i}-\frac{2 \sigma+\varepsilon r_{3}^{2}}{R^{2}}x^{i},
x^{i}-\frac{2 \sigma+\varepsilon
r_{3}^{2}}{R^{2}}y^{i},\frac{\varepsilon
r^{2}-2\sigma}{R^{2}}z^{j}) \leqno(3.11)\] for $i \in \{1,...,p\}$
and $j \in \{1,...,q\}$.

If we decompose $\PP X$  in the tangential and normal components
at sphere $S^{2p+q-1}$ (where X is a tangent vector field at
$S^{2p+q-1}(R)$), we obtain
\[ \PP X = \overline{P}X+\overline{u}(X)N_{1} \leqno(3.12) \]

From $\overline{u}(X)=<X,\overline{\xi}>$ and (3.7) we obtain
\[\overline{u}(X)=\frac{1}{R}[\sum_{i=1}^{p}(x^{i}Y^{i}+y^{i}X^{i})
+\varepsilon \sum_{j=1}^{q}z^{j}Z^{j}]\leqno(3.13)\] From $\PB
X=\PP X -\overline{u}(X)N_{1}$ we obtain
\[\: \PB X=(Y^{i}-\frac{\overline{u}(X)}{R}x^{i},
X^{i}-\frac{\overline{u}(X)}{R}y^{i},\varepsilon
Z^{j}-\frac{\overline{u}(X)}{R}z^{j}) \leqno(3.14) \] where
$X:=(X^{i},Y^{i},Z^{j})$ is a tangent vector at sphere in any
point  $(x^{i},y^{i},z^{j})$ and $\overline{u}(X)$ was defined in
(3.13).

Therefore, from the relations (3.10), (3.11), (3.13), (3.14) we
have a $(\PB,<>,\overline{\xi},\overline{u},a_{11})$ induced
structure by $\PP$ from $E^{2p+q}$ on the sphere $S^{2p+q-1}(R)$
of codimension 1 in Euclidean space $E^{2p+q}$. We use the
following notation $\overline{a}:=a_{11}$. If $a_{11}\neq 1$ in
every point of $S^{2p+q-1}(R)$, from the Theorem $2.2$ we obtain
that the structure
$(\PB,<>,\overline{\xi},\overline{u},\overline{a})$ is an
$(\overline{a},1)f$ normal structure, because $S^{2p+q-1}(R)$ is a
totally umbilical hypersurface in $E^{2p+q}$ and from this we have
that the tensor field P commutes with the Weingarten operator A.

%>>>>>>>>>>>>>>>>>>>>>>>>>>>>>>>>>>>>>>>>>>>>>>>>>>>>>>>>>>>>>>>>>>>>>>>>>>>>>
\textbf{ Example 2.} In this example, we construct an
$(\widehat{a},1)f$-structure on the product of spheres
$S^{2p-1}(r) \times S^{q-1}(r_{3})$.

Let $E^{2p+q}$ ($p,q \in \mathbb{N}^{*}$) be Euclidean space ($p,q
\in \mathbb{N}^{*}$) endowed with an almost product Riemannian
structure $\PP$ defined in $(3.1)$. It is obvious that
$E^{2p+q}=E^{2p} \times E^{q}$ and in each of spaces $E^{2p}$ and
$E^{q}$ respectively, we can get a hypersphere
$S^{2p-1}(r)=\{(x^{1},...,x^{p},y^{1},...,y^{p}),\sum_{i=1}^{p}((x^{i})^{2}+(y^{i})^{2})=r^{2}\}$
and
$S^{q-1}(r_{3})=\{(z^{1},...,z^{q}),\sum_{j=1}^{q}(z^{j})^{2}=r_{3}^{2}\}$
respectively, where $r^{2}+r_{3}^{2}=R^{2}$. We construct (in the
same manner like in \cite{IS} or \cite{eu7}) the product manifold
$S^{2p-1}(r) \times S^{q-1}(r_{3})$. Any point of $S^{2p-1}(r)
\times S^{q-1}(r_{3})$ has the coordinates
$(x^{1},...,x^{p},y^{1},...,y^{p},z^{1},...,z^{q}):=(x^{i},y^{i},z^{j})$
which verifies (3.3). Thus, $S^{2p-1}(r) \times S^{q-1}(r_{3})$ is
a submanifold of codimension 2 in $E^{2p+q}$. Furthermore,
$S^{2p-1}(r) \times S^{q-1}(r_{3})$ is a submanifold of
codimension 1 in $S^{2p+q-1}(R)$. Therefore we have the successive
imbedded
\[S^{2p-1}(r) \times S^{q-1}(r_{3})
 \hookrightarrow S^{2p+q-1}(R) \hookrightarrow E^{2p+q} \]
The tangent space in a point $(x^{i},y^{i},z^{j})$ at the product
of spheres $S^{2p-1}(r) \times S^{q-1}(r_{3})$ is
$T_{(x^{1},...,x^{p},y^{1},...,y^{p},\underbrace{o,...,o}_{q})}
S^{2p-1}(r)
 \oplus
 T_{(\underbrace{o,...,o}_{2p},z^{1},...,z^{q})}S^{q-1}(r_{3})$.

A vector $(X^{1},...,X^{p},Y^{1},...,Y^{p})$ from
$T_{(x^{1},...,x^{p},y^{1},...,y^{p})}E^{2p}$ is tangent to
$S^{2p-1}(r)$ if and only if
\[\sum_{i=1}^{p}x^{i}X^{i}+\sum_{i=1}^{p}y^{i}Y^{i}=0 \leqno(3.15)\]
and it can be identified by
$(X^{1},...,X^{p},Y^{1},...,Y^{p},\underbrace{0,...,0}_{q})$ from
$E^{2p+q}$. A vector $(Z^{1},...,Z^{q})$ from
$T_{(z^{1},...,z^{q})}E^{q}$ is tangent to $S^{q-1}(r_{3})$ if and
only if
\[\sum_{j=1}^{q}z^{j}Z^{j}=0 \leqno(3.16)\]
and it can be identified by
$(\underbrace{0,...,0}_{2p},Z^{1},...,Z^{q})$ from $E^{2p+q}$.
Consequently, for any point $(x^{i},y^{i},z^{j}) \in S^{2p-1}(r)
\times S^{q-1}(r_{3})$ we have $(X^{i},Y^{i},Z^{j}) \in
T_{(x^{1},...,x^{p},y^{1},...,y^{p},z^{1},...,z^{q})}(S^{2p-1}(r)
\times S^{q-1}(r_{3}))$ if and only if the relations (3.15) and
(3.16) are satisfied. Furthermore, we remark that
$(X^{i},Y^{i},Z^{j})$ is a tangent vector field at $S^{2p+q-1}(R)$
and from this it follows that
\[T_{(x^{i},y^{i},z^{j})}(S^{2p-1}(r) \times S^{q-1}(r_{3})) \subset
T_{(x^{i},y^{i},z^{j})}S^{2p+q-1}(R),\] for any point
$(x^{i},y^{i},z^{j}) \in S^{2p-1}(r) \times S^{q-1}(r_{3})$.

The normal unit vector field $N_{1}$ at $S^{2p+q-1}(R)$ is also a
normal vector field at $(S^{2p-1}(r) \times S^{q-1}(r_{3}))$, when
it is considered in its points. We construct an unit vector field
$N_{2}$ on $S^{2p+q-1}$, by the relation
\[N_{2}=\frac{1}{R}(\frac{r_{3}}{r}x^{i},\frac{r_{3}}{r}y^{i},-\frac{r}{r_{3}}z^{j})
\leqno(3.17)\] We remark that $N_{2}$ is orthogonal on $N_{1}$.
From (3.15) and (3.16) we obtain that $N_{2}$ is orthogonal at
$S^{2p-1}(r) \times S^{q-1}(r_{3})$. Thus, $(N_{1},N_{2})$ is an
orthonormal basis in $T_{(x^{i},y^{i},z^{j})}^{\bot}S^{2p-1}(r)
\times S^{q-1}(r_{3})$ in any point $(x^{i},y^{i},z^{j}) \in
S^{2p-1}(r) \times S^{q-1}(r_{3})$.

From the decomposed of $\PP N_{k}$ in tangential and normal
components at $S^{2p-1}(r) \times S^{q-1}(r_{3})$, we obtain
\[\PP N_{\A}=\widehat{\xi_{\A}}+a_{\A 1}N_{1}+a_{\A 2}N_{2}, \:\: \A \in \{1,2\}. \leqno(3.18)\]
$a_{11}$ was defined in (3.10). From $(2.3)(iii)$ we have
$a_{12}=a_{21}$ and from $a_{\A\B}=<\PP N_{\A},N_{\B}>$ ($\A,\B
\in \{1,2\}$), we obtain
\[a_{12}=a_{21}=\frac{(2 \sigma - \varepsilon r^{2})r_{3}}{rR^{2}} \leqno(3.19)\]
and
\[a_{22}=\frac{2 \sigma r_{3}^{2}+ \varepsilon r^{4}}{r^{2}R^{2}} \leqno(3.20)\]
From (3.18), (3.19), (3.20) and (3.10) we obtain
\[\widehat{\xi_{1}}=\frac{1}{R}(y^{i}-\frac{2 \sigma}{r^{2}}x^{i},
x^{i}-\frac{2 \sigma}{r^{2}}y^{i},0). \leqno(3.21)\] and
\[\widehat{\xi_{2}}=\frac{r_{3}}{rR}(y^{i}-\frac{2 \sigma}{r^{2}}x^{i},
x^{i}-\frac{2 \sigma}{r^{2}}y^{i},0). \leqno(3.22)\]

Therefore, the matrix $\widehat{a}:=(a_{\A\B})_{2}$ is given by
\[\widehat{a}:= \begin{pmatrix}
   \frac{2 \sigma + \varepsilon r_{3}^{2}}{R^{2}} &  \frac{(2 \sigma - \varepsilon r^{2})r_{3}}{rR^{2}}\\
   \frac{(2 \sigma - \varepsilon r^{2})r_{3}}{rR^{2}} & \frac{2 \sigma r_{3}^{2}+ \varepsilon r^{4}}{r^{2}R^{2}}
\end{pmatrix} \leqno(3.23)\]

From the decomposed of $\PP X$ in tangential and normal components
at $S^{2p-1}(r) \times S^{q-1}(r_{3})$ (where
$X:=(X^{i},Y^{i},Z^{j})$ is a tangent vector field on $S^{2p-1}(r)
\times S^{q-1}(r_{3})$), we obtain
\[\PP X = \widehat{P} X+\widehat{u}_{1}(X)N_{1}+\widehat{u}_{2}(X)N_{2}. \leqno(3.24) \]
We use the following notation
\[\tau:=\sum_{i=1}^{p}(x^{i}Y^{i}+y^{i}X^{i}) \leqno(3.25)\] From
$\widehat{u}_{\A}(X)=<X,\widehat{\xi}_{\A}>$ (with $\A \in
\{1,2\}$) and (3.15) we obtain
\[\widehat{u}_{1}(X)=\frac{\tau}{R},\quad
\widehat{u}_{2}(X)=\frac{r_{3}\tau}{rR} \leqno(3.26)\] for any
tangent vector $X$ on $S^{2p-1}(r) \times S^{q-1}(r_{3})$ in a
point $(x^{i},y^{i},z^{i}) \in S^{2p-1}(r) \times S^{q-1}(r_{3})$.
From (3.24) and (3.26), we obtain the tangent component
$\widehat{P}X$ of $\PP X$ on $S^{2p-1}(r) \times S^{q-1}(r_{3})$,
which is given by
\[\widehat{P} X
=(Y^{i}-\frac{\tau}{r^{2}}x^{i},
X^{i}-\frac{\tau}{r^{2}}y^{i},\varepsilon Z^{j}) \leqno(3.27)\]
Consequently, we obtain the
$(\widehat{P},<>,\widehat{\xi}_{1},\widehat{\xi}_{2},\widehat{u}_{1},\widehat{u}_{2},\widehat{a})$
induced structure on the product of spheres $S^{2p-1}(r) \times
S^{q-1}(r_{3})$ by the almost product Riemannian structure
$(\PP,<,>)$ which is an $(\widehat{a},1)f$ Riemannian structure
induced on the submanifold $S^{2p-1}(r) \times S^{q-1}(r_{3})$ of
codimension 2 in Euclidean space $E^{2p+q}$.

%>>>>>>>>>>>>>>>>>>>>>>>>>>>>>>>>>>>>>>>>>>>>>>>>>>>>>>>>>>>>>>>>>>>>>>>>>>>>>
\textbf{ Example 3.} In this example, we construct an
$(a,1)f$-structure on the product of spheres $S^{p-1}(r_{1})
\times S^{p-1}(r_{2})\times S^{q-1}(r_{3})$.

Let $E^{2p+q}$ be the Euclidean space endowed with an almost
product Riemannian structure $\PP$ which was defined in (3.1). We
have $E^{2p+q}=E^{p} \times E^{p} \times E^{q}$ and in each of
spaces $E^{p}$ we can get a hypersphere
\[S^{p-1}(r_{1})=\{(x^{1},...,x^{p}), \: \sum_{i=1}^{p}(x^{i})^{2}=r_{1}^{2}\}, \leqno(3.28)\]
and
\[S^{p-1}(r_{2})=\{(y^{1},...,y^{p}), \: \sum_{i=1}^{p}(y^{i})^{2}=r_{2}^{2}\}, \leqno(3.29)\]
respectively, with the assumption that
$r_{1}^{2}+r_{2}^{2}=r^{2}$.

We construct the product manifold  $S^{p-1}(r_{1}) \times
S^{p-1}(r_{2})\times S^{q-1}(r_{3})$.

Any point x of $S^{p-1}(r_{1}) \times S^{p-1}(r_{2}) \times
S^{q-1}(r_{3})$ has the coordinates
\[(x^{1},...,x^{p},y^{1},...,y^{p},z^{1},...,z^{q}):=(x^{i},y^{i},z^{j})\]
which verifies (3.3). Thus, $S^{p-1}(r_{1})\times S^{p-1}(r_{2})
\times S^{q-1}(r_{3})$ is a submanifold of codimension 3 in
$E^{2p+q}$. Furthermore, $S^{p-1}(r_{1})\times S^{p-1}(r_{2})
\times S^{q-1}(r_{3})$ is a submanifold of codimension $1$ in
$S^{2p-1}(r) \times S^{q-1}(r_{3})$ and a submanifold of
codimension $2$ in $S^{2p+q-1}(R)$. Therefore, we have the
successive imbedded
\\ $S^{p-1}(r_{1})\times S^{p-1}(r_{2}) \times
S^{q-1}(r_{3})\hookrightarrow
 S^{2p-1}(r) \times S^{q-1}(r_{3})
 \hookrightarrow S^{2p+q-1}(R) \hookrightarrow E^{2p+q}$

The tangent space in a point $x:=(x^{i},y^{i},z^{j})$ at the
product of spheres $S^{p-1}(r_{1})\times S^{p-1}(r_{2}) \times
S^{q-1}(r_{3})$ is
\[T_{(x^{i},0,0)} S^{p-1}(r_{1}) \oplus T_{(0,y^{i},0)} S^{p-1}(r_{2})
 \oplus T_{(0,0,z^{j})}S^{q-1}(r_{3})\]
where $(x^{i},0,0)$ is a notation for
$(x^{1},...,x^{p},\underbrace{0,...,0}_{p+q})$, $(0,y^{i},0)$ is a
notation for
$(\underbrace{0,...,0}_{p},y^{1},...,y^{p},\underbrace{0,...,0}_{q})$
and $(0,0,z^{j})$ is a notation for
$(\underbrace{0,...,0}_{2p},z^{1},...,z^{q})$.

A vector $(X^{1},...,X^{p})$ from $T_{(x^{1},...,x^{p})}E^{p}$ is
tangent to $S^{p-1}(r_{1})$ if and only if
\[\sum_{i=1}^{p}x^{i}X^{i}=0 \leqno(3.30)\]
and it can be identified by
$(X^{1},...,X^{p},\underbrace{0,...,0}_{p+q})$ from $E^{2p+q}$.

A vector $(Y^{1},...,Y^{p})$ from $T_{(y^{1},...,y^{p})}E^{p}$ is
tangent to $S^{p-1}(r_{2})$ if and only if
\[\sum_{i=1}^{p}y^{i}Y^{i}=0 \leqno(3.31)\]
and it can be identified by
$(\underbrace{0,...,0}_{p},Y^{1},...,Y^{p},\underbrace{0,...,0}_{q})$
from $E^{2p+q}$.

Consequently, for any point $(x^{i},y^{i},z^{j}) \in
S^{p-1}(r_{1}) \times S^{p-1}(r_{2}) \times S^{q-1}(r_{3})$ we
have
\[(X^{i},Y^{i},Z^{j}) \in
T_{(x^{i},y^{i},z^{j})}(S^{p-1}(r_{1}) \times S^{p-1}(r_{2})
\times S^{q-1}(r_{3}))\] if and only if the relations (3.16),
(3.30) and (3.31) are satisfied. Furthermore, we remark that
$(X^{i},Y^{i},Z^{j})$ is a tangent vector field at the product of
spheres $S^{2p-1}(r) \times S^{q-1}(r_{3})$ and at
$S^{2p+q-1}(R)$, respectively. Thus, it follows that \\
$T_{(x^{i},y^{i},z^{j})}(S^{p-1}(r_{1}) \times S^{p-1}(r_{2})
\times S^{q-1}(r_{3})) \subset T_{(x^{i},y^{i},z^{j})}(S^{2p-1}(r)
\times S^{q-1}(r_{3}))$, for any point $(x^{i},y^{i},z^{j}) \in
S^{p-1}(r_{1}) \times S^{p-1}(r_{2}) \times S^{q-1}(r_{3})$.

The normal unit vector fields $N_{1}$ and $N_{2}$ at $(S^{2p-1}(r)
\times S^{q-1}(r_{3}))$ are also normal vector fields at
$(S^{p-1}(r_{1}) \times (S^{p-1}(r_{2}) \times S^{q-1}(r_{3}))$,
when they are considered in its points. We construct a unit vector
field $N_{3}$ by the relation
\[N_{3}=\frac{1}{r}(\frac{r_{2}}{r_{1}}x^{i},-\frac{r_{1}}{r_{2}}y^{i},0)
\leqno(3.32)\] We remark that $N_{3}$ is orthogonal on $N_{1}$ and
$N_{2}$, respectively. Furthermore, from (3.16), (3.30) and (3.31)
we obtain that $N_{3}$ is orthogonal at $S^{p-1}(r_{1}) \times
S^{p-1}(r_{1}) \times S^{q-1}(r_{3})$. Thus, $(N_{1},N_{2},
N_{3})$ is an orthonormal basis in $T_{(x^{i},y^{i},z^{j})}^{\bot}
S^{p-1}(r_{1}) \times S^{p-1}(r_{2}) \times S^{q-1}(r_{3})$ in any
point $(x^{i},y^{i},z^{j}) \in S^{p-1}(r_{1}) \times
S^{p-1}(r_{2}) \times S^{q-1}(r_{3})$.

From the decomposed of $\PP N_{\A}$ (with $\A \in \{1,2,3\}$) in
the tangential and normal components at $S^{p-1}(r_{1}) \times
S^{p-1}(r_{2}) \times S^{q-1}(r_{3})$, we obtain
\[\PP N_{\A}=\xi_{\A}+a_{\A 1}N_{1}+a_{\A 2}N_{2}+a_{\A 3}N_{3}, \leqno(3.33)\]
for every $ \A \in \{1,2,3\}$.  From (3.10) we have $a_{11}$, from
(3.19) we have $a_{12}=a_{21}$ and from (3.20) we have $a_{22}$.
Furthermore, from $(2.3)(iii)$ we have that $a_{13}=a_{31}$ and
$a_{23}=a_{32}$.

It is obvious that \[a_{\A\B}=<\PP N_{\A},N_{\B}>, \quad (\A,\B
\in \{1,2,3\}). \leqno(3.34)\]

From this, we have
\[a_{13}=a_{31}=\frac{r_{2}^{2}-r_{1}^{2}}{r_{1}r_{2}rR}\sigma \leqno(3.35)\]
\[a_{23}=a_{32}=\frac{(r_{2}^{2}-r_{1}^{2})r_{3}}{r_{1}r_{2}r^{2}R}\sigma \leqno(3.36)\]
and
\[a_{33}=-\frac{2\sigma}{r^{2}} \leqno(3.37)\]

From (3.33) we obtain
\[\xi_{1}=\frac{1}{R}(y^{i}-\frac{\sigma}{r_{1}^{2}}x^{i},
x^{i}-\frac{\sigma}{r_{2}^{2}}y^{i},0). \leqno(3.38)\]
\[\xi_{2}=\frac{r_{3}}{rR}(y^{i}-\frac{\sigma}{r_{1}^{2}}x^{i},
x^{i}-\frac{\sigma}{r_{2}^{2}}y^{i},0). \leqno(3.39)\] and
\[\xi_{3}=\frac{1}{r}(\frac{\sigma}{r_{1}r_{2}}x^{i}-\frac{r_{1}}{r_{2}}y^{i},
\frac{r_{2}}{r_{1}}x^{i}-\frac{\sigma}{r_{1}r_{2}}y^{i},0).
\leqno(3.40)\]

Therefore, the matrix $a:=(a_{\A\B})_{3}$
 is given by
\[a =\begin{pmatrix}
    \frac{2 \sigma + \varepsilon r_{3}^{2}}{R^{2}} &
    \frac{(2 \sigma - \varepsilon r^{2})r_{3}}{rR^{2}} &
    \frac{r_{2}^{2}-r_{1}^{2}}{r_{1}r_{2}rR}\sigma \\
   \frac{(2 \sigma - \varepsilon r^{2})r_{3}}{rR^{2}} &
   \frac{2 \sigma r_{3}^{2}+ \varepsilon r^{4}}{r^{2}R^{2}}&
   \frac{(r_{2}^{2}-r_{1}^{2})r_{3}}{r_{1}r_{2}r^{2}R}\sigma \\
   \frac{r_{2}^{2}-r_{1}^{2}}{r_{1}r_{2}rR}\sigma &
    \frac{(r_{2}^{2}-r_{1}^{2})r_{3}}{r_{1}r_{2}r^{2}R}\sigma &
    -\frac{2\sigma}{r^{2}}
\end{pmatrix} \leqno(3.41)\]
From the decomposed of $\PP X$ in tangential and normal components
at the product of spheres $S^{p-1}(r_{1}) \times S^{p-1}(r_{2})
\times S^{q-1}(r_{3})$ (where $X:=(X^{i},Y^{i},Z^{j})$ is a
tangent vector field on $S^{p-1}(r_{1}) \times S^{p-1}(r_{2})
\times S^{q-1}(r_{3})$), we obtain
\[\PP X = P X+u_{1}(X)N_{1}+u_{2}(X)N_{2}+u_{3}(X)N_{3}. \leqno(3.42) \]
From $(2.3)(vi)$ we have $u_{\A}(X)=<X,\xi_{\A}>$ (with $\A \in
\{1,2,3\}$) and using (3.37) and (3.38) we obtain
\[u_{1}(X)=\frac{\tau}{R}, \leqno(3.43)\]
\[u_{2}(X)=\frac{r_{3}\tau}{rR} \leqno(3.44)\] and
respectively
\[u_{3}(X)=\frac{1}{r}[\frac{r_{2}}{r_{1}}\sum_{i=1}^{p}x^{i}Y^{i}-
\frac{r_{1}}{r_{2}}\sum_{i=1}^{p}y^{i}X^{i}]. \leqno(3.45)\] for
any tangent vector $X:=(X^{i},Y^{i},Z^{j})$  in any point
$(x^{i},y^{i},z^{j})$ on the product of spheres  $S^{p-1}(r_{1})
\times S^{p-1}(r_{2}) \times S^{q-1}(r_{3})$, with the notation
$(3.25)$ for $\tau$. From (3.42), (3.43), (3.44) and (3.45), we
obtain the tangent component $P X$ of $\PP X$ in any point
$(x^{i},y^{i},z^{j})$ on $S^{p-1}(r_{1}) \times S^{p-1}(r_{2})
\times S^{q-1}(r_{3})$, which has the form
\[P X =(Y^{i}-\frac{\sum_{k=1}^{p}x^{k}Y^{k}}{r_{1}^{2}}x^{i},
X^{i}-\frac{\sum_{k=1}^{p}y^{k}X^{k}}{r_{2}^{2}}y^{i},\varepsilon
Z^{j}) \leqno(3.46)\] Consequently, we obtain the induced
structure
$(P,<>,\xi_{1},\xi_{2},\xi_{3},u_{1},u_{2},u_{3},(a_{\A\B})_{3})$
on $S^{p-1}(r_{1}) \times S^{p-1}(r_{2}) \times S^{q-1}(r_{3})$ by
the almost product Riemannian structure $(\PP,<>)$, which is an
$(a,1)f$ Riemannian structure induced on the product of spheres
$S^{p-1}(r_{1}) \times S^{p-1}(r_{2}) \times S^{q-1}(r_{3})$ in
$E^{2p+q}$.

\rem If we consider the almost product structure on $E^{2p+q}$
given by $\PP:E^{2p+q}\rightarrow E^{2p+q}$,
\[ \PP(x^{1},...,x^{p},y^{1},...,y^{p},z^{1},...,z^{q})= \]
\[=(y^{1},...,y^{p},x^{1},...,x^{p},\varepsilon_{1}
z^{1},...,\varepsilon_{q} z^{q})\] where $\varepsilon_{j} = \pm 1$
(for every $j\in \{1,...,q\})$, we obtain a generalization of
$(a,1)f$ induced structures defined in examples from above.

%>>>>>>>>>>>>>>>>>>>>>>>>>>>>>>>>>>>>>>>>>>>>>>>>>>>>>>>>>>>>>>>>>>>>>>>>>>>>>>>>>>>>

\bigskip

 "\c{S}tefan cel Mare" University, Suceava, Romania\\
  E-mail address: cristinah@usv.ro

\end{document}